\pdfoutput=1
\documentclass{article}
\usepackage[margin=25mm]{geometry}
\usepackage{amsmath}
\usepackage{amsfonts}
\usepackage{amssymb}
\usepackage{graphicx}
\usepackage[utf8]{inputenc}
\usepackage[T2A]{fontenc}
\usepackage[ukrainian,russian,english]{babel}
\usepackage[pdftex]{hyperref}
\usepackage{xcolor}
\usepackage[all]{xy}
\newtheorem{lemma}{Lemma}
\newtheorem{theorem}{Theorem}

\title{
Structures of the flows with a unique singular point on the 2-dimensional disk }

\author{Alexandr Prishlyak and Serhii Stas}

\begin{document}

\maketitle
\begin{abstract}
   We investigate topological propeties of flows with one singular point and without closed orbits on the 2-dimensional disk.  To classify such flows, destingueshed graph is used, which is a two-colored rooted tree imbedded in the plane.   We construct a code of the flow and have found all possible structures of the flows with no more then 7 sepapratrices. 
\end{abstract}
\section*{Introduction}

On closed manifolds, if we consider Morse functions that are functions in general position, then the number of critical points of each index satisfies the Morse inequalities. For example, it is not less than the minimum number of generators of the corresponding homology group. If a function has two critical points, then the manifold is a sphere. On closed surfaces, there are functions with three critical points (one of which is degenerate). If we consider the corresponding gradient flows, then we have the corresponding restrictions on the number of singular points. There are no such restrictions for arbitrary flows. In particular, there exists a flow with one singular point whose Poincare rotation index is equal to the Euler characteristic of the manifold. If the flow has a finite number of closed trajectories, then by introducing pairs of singular points of neighboring Morse indices on them, one can obtain a flow without singular points. The purpose of this paper is to investigate the topological structure of flows on a two-dimensional disk with one singular point and without closed orbits. It is obvious that this point have to lie on the boundary of the disk, since otherwise the boundary is a closed orbit.

The use of graph theory to classify dynamical systems is typical in two and three dimensions. 





The topological classifications of flows were obtained on closed 2-manifolds in \cite{bilun2023gradient, Kybalko2018, Oshemkov1998, Peixoto1973, prishlyak1997graphs, prishlyak2020three, akchurin2022three, prishlyak2022topological, prishlyak2017morse,  kkp2013, prishlyak2021flows,  prishlyak2020topology,   prishlyak2019optimal, prishlyak2022Boy}, 
 and on 2-manifolds with the boundary in
\cite{bilun2023discrete, bilun2023typical, loseva2016topology, prishlyak2017morse, prishlyak2022topological, prishlyak2003sum, prishlyak2003topological, prishlyak1997graphs, prishlyak2019optimal}.
Complete topological invariants of Morse-Smale flows on 3-manifolds was constructed in \cite{prish1998vek,  prish2001top, Prishlyak2002beh2, prishlyak2002ms,  prishlyak2007complete, hatamian2020heegaard, bilun2022morse, bilun2022visualization}.

Morse flows are gradient flows of Morse functions in general position.  The flow determinate the topological structure of the function if one fixes the value of functions in singular points\cite{lychak2009morse, Smale1961}. Therefore, Morse--Smale flows classification is related to the classification of the Morse functions.

Topological invariants of functions on oriented 2-maniofolds were constructed in \cite{Kronrod1950} and \cite{Reeb1946} and  in \cite{lychak2009morse} for  non-orientable two-dimensional manifolds, in   \cite{Bolsinov2004, hladysh2017topology, hladysh2019simple} for manifolds with boundary, in \cite{prishlyak2002morse} for non-compact manifolds. 

Topological invariants of smooth function on closed 2-manifolds was also investigated in \cite{bilun2023morseRP2, hladysh2019simple, hladysh2017topology,  prishlyak2002morse, prishlyak2000conjugacy,  prishlyak2007classification, lychak2009morse, prishlyak2002ms, prish2015top, prish1998sopr,  bilun2002closed,  Sharko1993}, on 2-manifolds with the boundary in \cite{hladysh2016functions, hladysh2019simple, hladysh2020deformations} and on closed 3- and 4-manifolds in  \cite{prishlyak1999equivalence, prishlyak2001conjugacy}.



We are looking at flows on the 2-disk $D^2=\{(x,y)\in \mathbb{R}^2 : x^2+y^2 \le 1\}$.

We describe all possible topological structures of flows with a  unique singular point on the boundary  and small number of separatrices.
We call such flows as 1-flows.

In the first section, we give the basic definitions and topological properties 0f 1-flows. In the second section, we describe the properties of separatrix diagrams and distinguishing graph, with is dual to the sets of separatrix. 
In the third section we introduce the code to classify the 1-flows.

\section{Topological properties of 1-flows}

Let us describe the topological properties of 1-flows.



\begin{figure}[h]
\center{\includegraphics[width=0.8\linewidth]{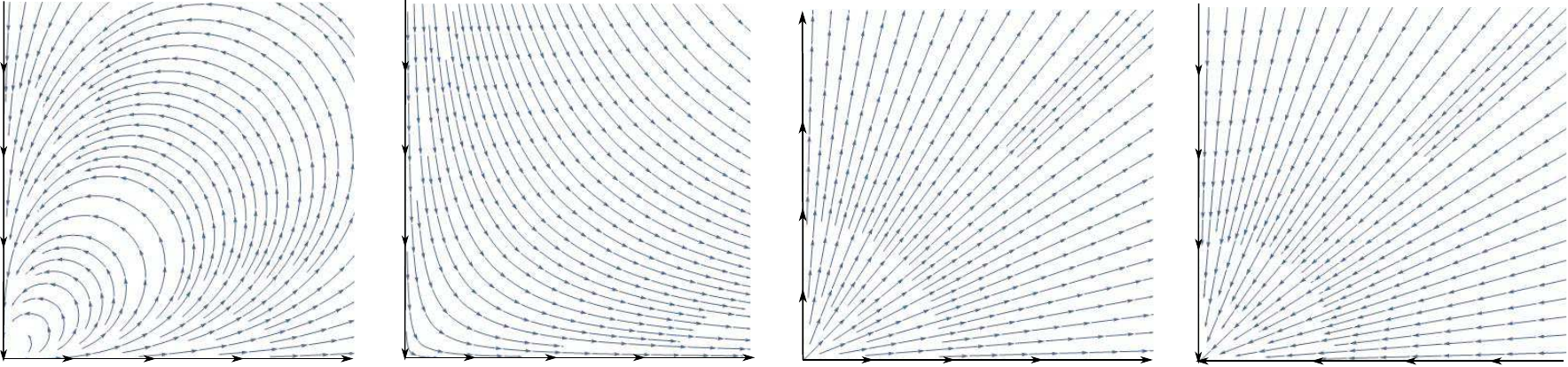}}
\center{1) \ \ \ \ \ \ \ \ \ \ \ \ \ \ \ \ \ \ \ \ \ \ \ \ \ 
2) \ \ \ \ \ \ \ \ \ \ \ \ \ \ \ \ \ \ \ \ \ \ \ \ \ 
3) \ \ \ \ \ \ \ \ \ \ \ \ \ \ \ \ \ \ \ \ \ \ \ \ \  4) }
\caption{types of vertex angles: 1) elliptic, 2) hyperbolic, 3) source, 4) sink}
\label{fixedpoint}
\end{figure}

Separatrices divide the 2-disk into regions, which we will call cells. Each cell is a curved polygon, the vertices of which are the singular point. Since there are no singular points inside the cell, each trajectory begins and ends at the singular point. There are 4 types of vertex angles: 1) elliptic, 2) hyperbolic, 3) sources, 4) sinks. In the neighbourhood of the elliptic point, the field is topologically equivalent to the field $\{x^3-3xy^2, 3x^2y-y^3\}, x\geq 0, y\geq 0$, in the vicinity of the hyperbolic point it is equivalent $\{ x,-y\}, x\geq 0, y\geq 0$, in the source to $\{ x,y\}, x\geq 0, y\geq 0$, in the sink to $\{ -x,-y\}, x\geq 0, y\geq 0$ see Fig. \ref{fixedpoint}.

All cells are of two types: 1) one of the corners is a source, another corner is a sink, and the rest are hyperbolic (polar cell); 2) one of the angles is elliptic, and the rest are hyperbolic (cyclic cell).
In a polar cell, trajectories begin and end at different points, and in a polar cell, at one. The boundary of the cyclic cell forms a cycle, and the polar one forms two oriented paths that start at the source and end at the sink.

\section{Separatrix diagram and  distinguishing  graph of 1-flow }

\begin{figure}[ht] 
\center{\includegraphics[width=0.85\linewidth]
{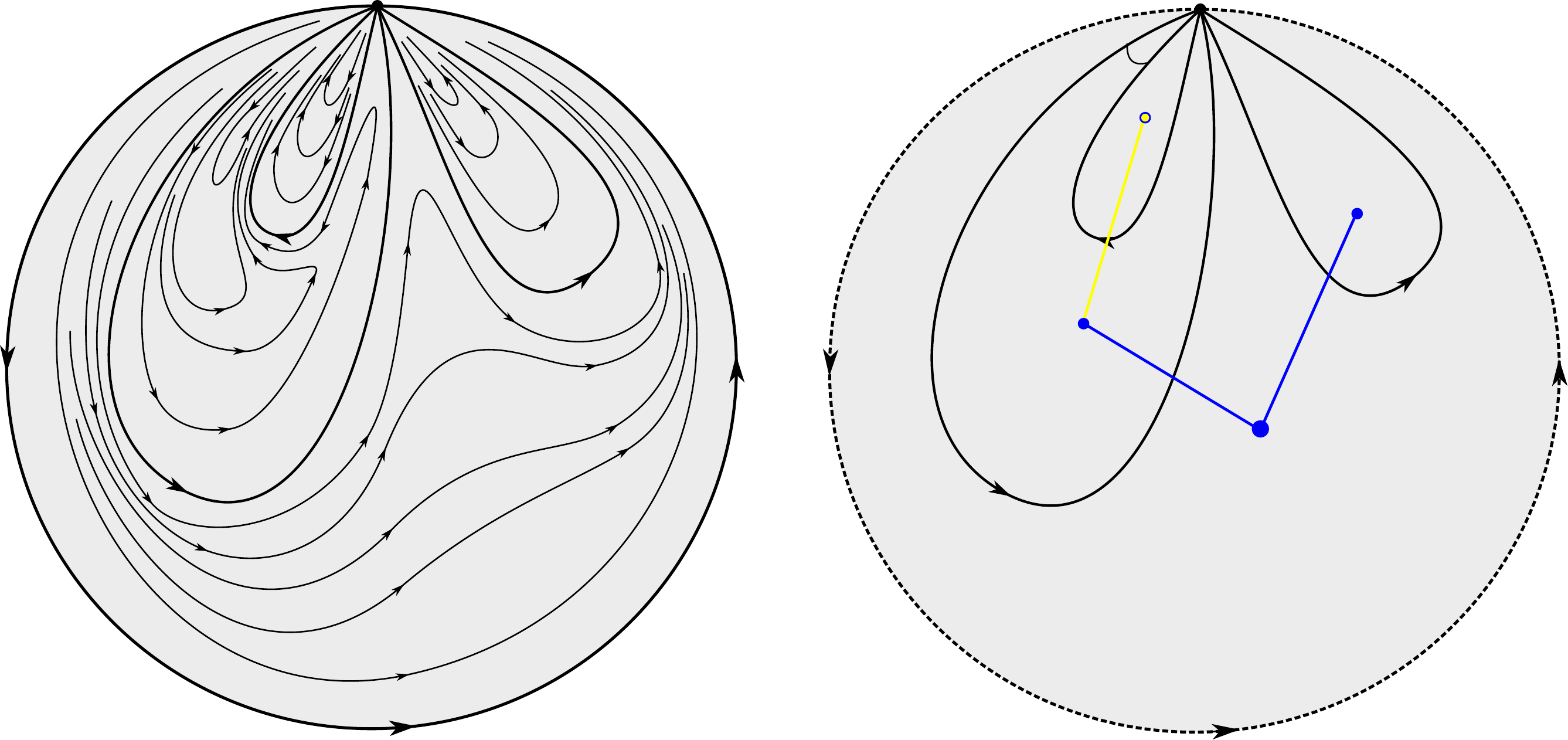}
}
\caption{A flow, separatrix diagram and destinguished graph with code $210\overline{0}$
}
\label{example}
\end{figure}

As is known, flows on a two-dimensional sphere can be topologically classified using a separatrix diagram. A 2-disk can be thought of as a region on a sphere, so we can use the separatrix diagram to classify such flows. In our case, the separatrix diagram consists of a graph embedded in a 2-disk, the only vertex of which is a singular point, and the edges are separatrices. In addition, each edge is oriented, and parabolic sectors are selected among the sectors in the vicinity of the singular point. An example of a 1-flow and its sepratrix diagram is shown in fig. \ref{example}.

Thus, a separatrix diagram of a 1-flow is a bunch of oriented circles with distinguished corners nested in a 2-disk.
For convenience, we replace the separatrix diagram with the dual graph. The highlighted vertex on the graph is the one that corresponds to the region adjacent to the boundary.
  
\begin{lemma} The graph dual to the separatrix diagram is a rooted tree embedded in the plane.
\end{lemma}
\textbf{Proof.} Let us consider the structure of a cell complex, in which the singular point is the 0-cell, and the remaining part of the boundary and the separatrix are 1-cells. The dual graph, as a 1-skeleton, defines the dual structure of a cellular complex on a 2-disk. Since the 1-skeleton of a surface with boundary is homotopy equivalent to it, then the dual graph is homotopy equivalent to a point, that is, a tree.

$ $

For each vertex, except for the root, there is one incident edge that belongs to a simple path from the root to the same vertex. Such an edge will be called the lower edge with respect to the given vertex. Other incident edges will be called upper. At the root, all incident edges are top. When drawing a rooted tree, the top edges will be shown higher from the top, and the bottom edge - below.

\begin{figure}[ht!] 
\center{\includegraphics[width=0.85\linewidth]{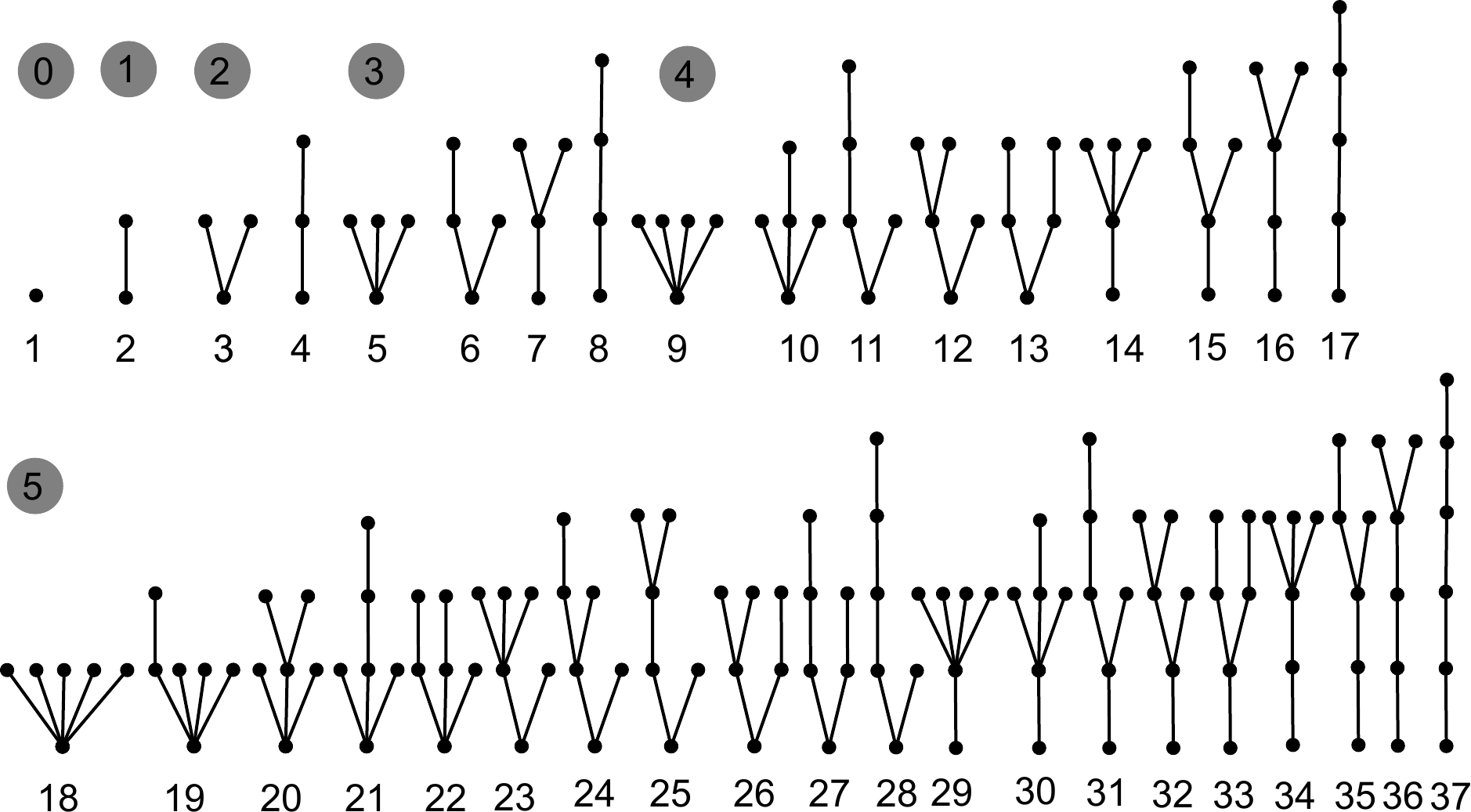}
}
\caption{rooted trees with at most five edges
}
\label{trees}
\end{figure}

For its bottom edge, the vertex is the top end, and for each of the top edges, the vertex is the bottom end.

  All possible rooted trees with no more than six vertices are shown in Fig. \ref{trees}.

\begin{figure}[ht!] 
\center{\includegraphics[width=0.50\linewidth]{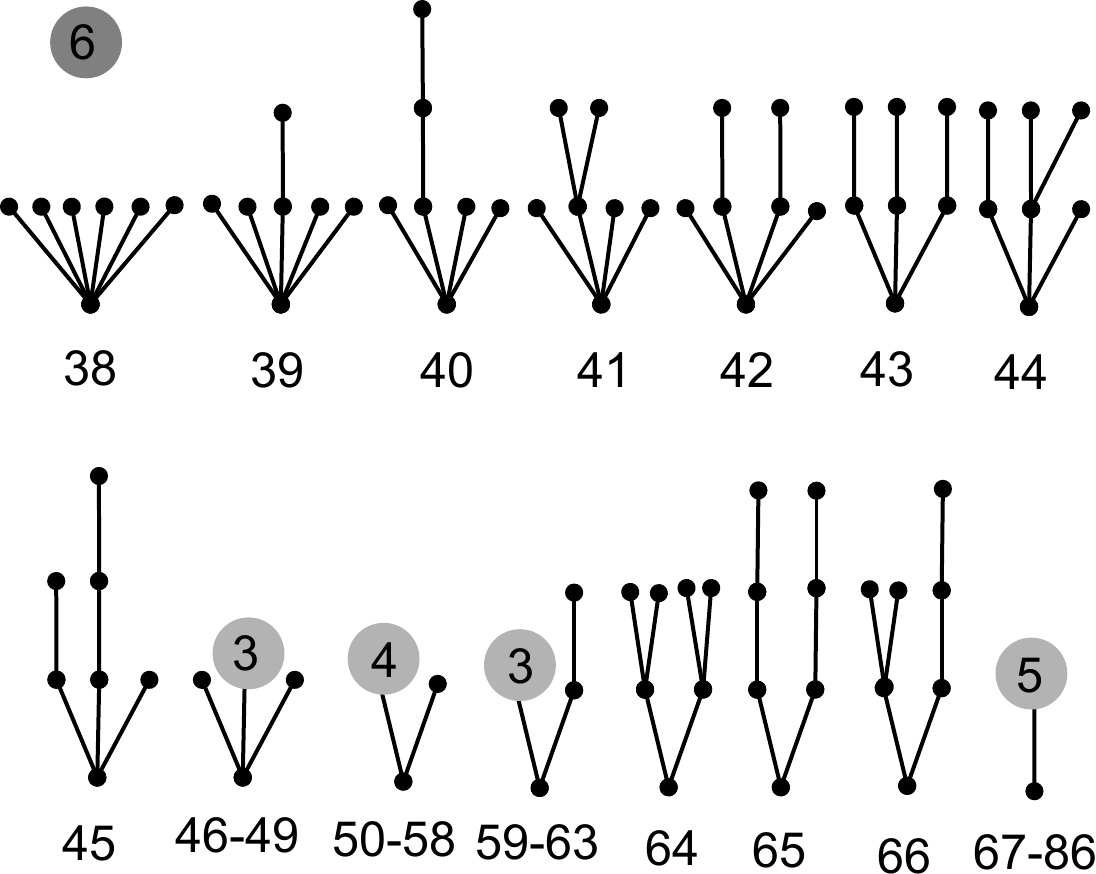}
}
\caption{rooted trees with six edges
}
\label{trees6}
\end{figure}

\begin{figure}[ht!] 
\center{\includegraphics[width=0.75\linewidth]{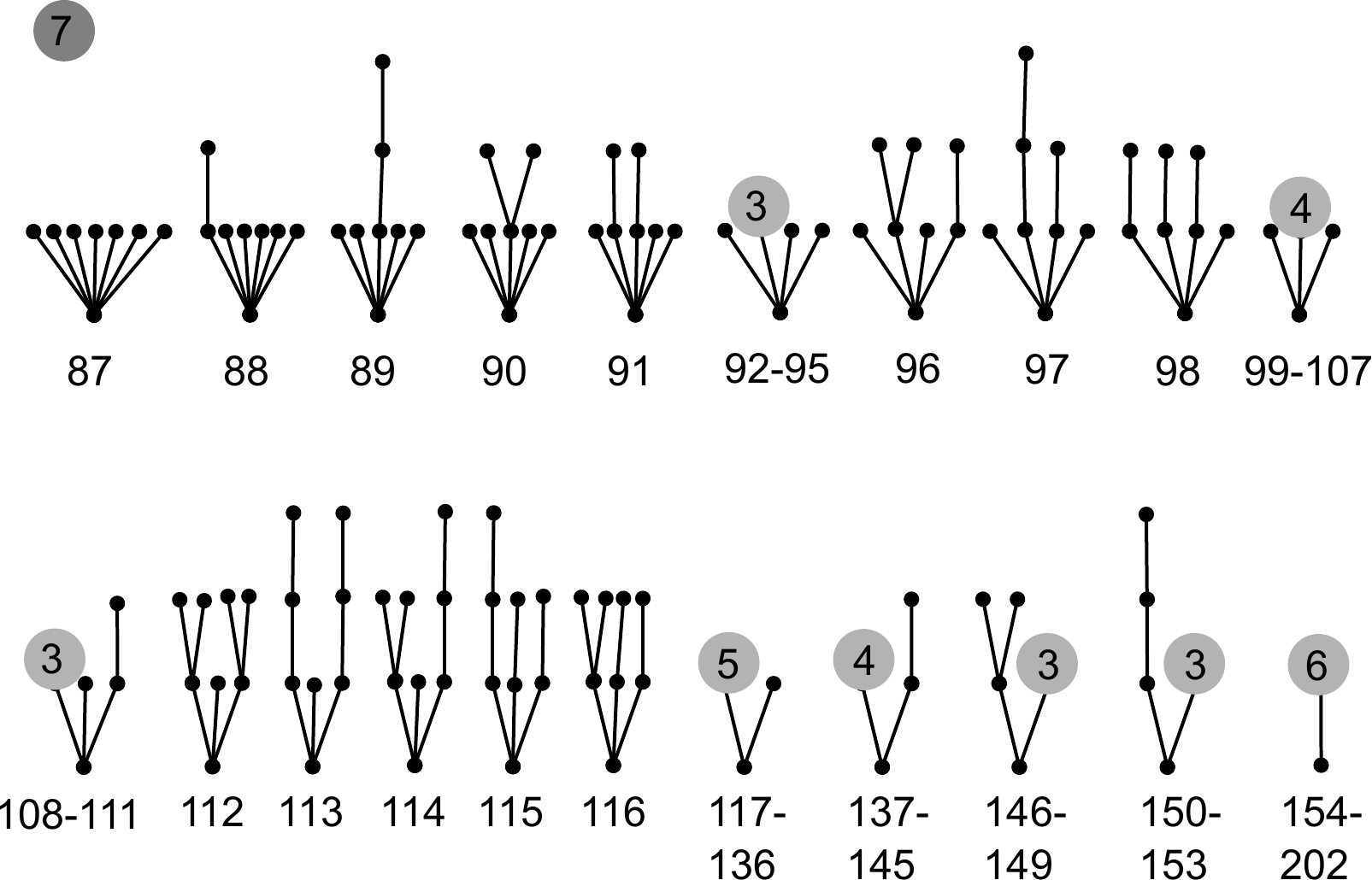}
}
\caption{rooted trees with seven edges
}
\label{trees7}
\end{figure}

Every tree has a canonical orientation of edges, given by the direction of movement from the lower to the upper end. Let us fix the orientation of the plane (two-dimensional disk) that induces the orientation of the flow at its boundary. For each edge of the tree, on the corresponding loop, there are two orientations - one coincides with the direction of the flow, and the second is generated by the orientation of the plane to the loop. If these orientations match, then color the edge black (1), and if they do not match, then color it red (-1).

If for some vertex all upper edges have the opposite color to the lower vertex, then the corresponding area will be cyclic. The flow in it will be set by specifying the elliptical angle. To do this, it is enough to indicate the loop that enters this corner. If the corresponding edge on the tree is upper, then we put a label on it, and if it is lower, then we leave it unchanged.
For each vertex, on the set of upper edges, we fix the order, in accordance with the traversal of the boundary of the corresponding region along the canonical orientation, starting from the lower loop.

A distinguishing flow graph is a two-color rooted tree with edge labels and orders of upper edges. Two distinguishing graphs are said to be equivalent if there is an isomorphism of them that preserves colors, labels and orders of upper edges and takes the root of the tree to the root.

\begin{theorem}
 Two 1-flows on a 2-disk are topologically equivalent if and only if their distinguishing graphs are equivalent.
\end{theorem}
\textbf{Proof.} Necessity follows from construction. Let us prove sufficiency. The orders on the set of upper edges define the rotation system on the graph. Then, in accordance with topological graph theory, there exists a unique, up to homeomorphism, embedding of a graph on a plane with this rotation system. Let's build a dual graph. It defines a set of separatrices. The orientation of the separatrices is given by the color of the corresponding edge in the tree, and the choice of the elliptical angle is given by the label.

\section{Code of 1-flow  }

To construct the code, we use a distinguishing 1-flow graph. The level of a vertex is the number of edges in the path between the vertex and the root. At level 0, there is one vertex, the root. Let's compose a sequence of integers equal to the numbers of the top edges of the vertices, starting from the vertex of level 0, then levels 1, 2, 3 $\ldots$ At each level, the vertices are ordered from left to right according to their location on the nested graph. If the lower edge of the vertex has the color -1, then put an upper line over the corresponding number, for example 3. If the lower edge has a label, then put a stroke after the corresponding number, for example. Thus, a 1-stream code is a set of integers, some of which may have an overbar or stroke, or a stroke and stroke at the same time. Note that the first number in the code has neither an overline nor a stroke.

We present codes for all possible flows with at most 3 separatrices:

1) 0;

2) 10, $1 \overline{0}$,  $1\overline{0}'$;

3) 200, $2\overline{0}0$, $20\overline{0}$,  
$2\overline{0}\overline{0}$, $2\overline{0}'\overline{0}$,  $2\overline{0}\overline{0}'$;

4) 110, $11\overline{0}$,  $11\overline{0}'$, 
$1\overline{1}\overline{0}$, $1\overline{1}'\overline{0}$, $1\overline{1}0$, $1\overline{1}'0$, 
$1\overline{1}0'$, $1\overline{1}'0'$.

For the trees shown in Fig. 3 (all edges have a canonical orientation and are not labeled) 1-flow codes are given in the second column of the table.

\begin{table}[ht!]

\begin{center}
\begin{tabular}{|c|c|c|c|c|}
\hline
tree & \ \ \ \ code \ \ \ \ & number of flows& tree embedings & total number\\
\hline
1 &0 & 1&1 &1 \\ \hline 
2 &10 & 3& 1&3 \\ \hline 
3 &200 & 6& 1& 6 \\ \hline 
4 & 110& 9& 1& 9\\ \hline 
5 & 3000& 10&1 &10 \\ \hline 
6 & 2100& 18& 2 & 36 \\ \hline 
7 & 1200&18 & 1&18 \\ \hline 
8 & 1110& 27& 1&27 \\ \hline 
9 & 40000& 15 & 1& 15 \\ \hline 
10 & 30100&30 & 3&90 \\ \hline 
11 & 21010& 54 & 2&108 \\ \hline 
12 & 22000& 36 & 2 & 72 \\ \hline 
13 & 21100& 54 & 1& 54 \\ \hline 
14 &13000& 30 & 1& 30\\ \hline 
15 &12100 & 54&2 &108 \\ \hline 
16 &11200 & 54& 1&54 \\ \hline 
17 &11110 & 81& 1&81\\ \hline 
18 &500000 &21 & 1&21 \\ \hline 
19 &410000 &45 &4 & 180 \\ \hline 
20 &302000 &60 &3 & 180\\ \hline
21 &301010 &90 &3 &270 \\ \hline 
22 &311000 &90 &3 &270 \\ \hline 
23 &230000 &60 &2 &120 \\ \hline 
24 &220100 &108 &4 &432 \\ \hline 
25 &210200 &108 &2 &216 \\ \hline 
26 &221000 &108 &2 &216 \\ \hline 
27 &211100 &162 &2 &324 \\ \hline 
28 &210110 &162 &2 &324 \\ \hline 
29 & 140000& 45 & 1& 45 \\ \hline 
30 & 130100&90 & 3&270 \\ \hline 
31 & 121010& 162 & 2&324 \\ \hline 
32 & 122000& 108 & 2 & 216 \\ \hline 
33 & 121100& 162 & 1& 162 \\ \hline 
34 &113000& 90 & 1& 90\\ \hline 
35 &112100 & 162&2 &324 \\ \hline 
36 &111200 & 162& 1&162 \\ \hline 
37 &111110 & 243& 1&243\\ \hline 

\end{tabular}
\end{center}
\caption{\label{tab1} number of 1-flow srtuctures}
\end{table} 

\begin{theorem}
Two 1-flows are topologically equivalent if and only if their codes are the same.
\end{theorem}
\textbf{Proof.} It follows from Theorem 1 that it suffices to prove the equivalence of distinguishing graphs. The necessity follows from the uniqueness of constructing the code according to the distinguishing tree. Let us prove sufficiency. To do this, it suffices to specify an algorithm for constructing a distinguishing graph from a code. We build the tree sequentially level by level from bottom to top. From the bottom root vertex, draw up as many edges as the first number of the code. Further, from the upper ends of the drawn edges from left to right, we draw up the edges in the amount indicated in the code. From the ends of the newly drawn edges, we draw new edges up as long as there are numbers in the code. The orientation and coloring of the edges (the orientation of the separatrices and the selection of the parabolic angle) are determined by the upper lines and strokes. Thus, the flow structure is unambiguously restored from the code.

$ $

For each, except the first, number in the code, there is a root, the number corresponding to the bottom cell, in which the corresponding edge enters as the top edge. We call two numbers connected if their edges are the top edges of the same cell. We call a number without overline and the same overline number by reverse numbers.

\begin{theorem}

An 1-flow code with n separatrices has properties:

1. The sum of all code numbers is equal to $n$. The number of all numbers is equal to $n+1$.

2. The first number does not have an overline and a prime.

3. For any $k$ not exceeding $n$, the sum of the first $k$ numbers is not less than $k$.

4. If some number has a prime, then all connected numbers do not have a prime, and all these numbers either have an overline, or all of them do not have it. Moreover, their root has or doesn't have overline opposite to the upper numbers.
\end{theorem}
\textbf{Proof. }
1. By definition, the sum of all numbers is equal to the number of upper edges. Since each edge is the top one for only one vertex, this sum is equal to the number of edges, i.e. $n$. Each code number corresponds to a vertex, so their total number is $n+1$.

2. The presence of an overline over the first number depends on the orientation of the area boundary. Since this orientation is always the same (counterclockwise), there is no overline above the first number. The prime near the first number is also missing, since the outer arc cannot be the upper arc in the elliptic region.

3. This property follows from the fact that the number of vertices that are located below or to the left of a given vertex is less than the number of edges that are located below this vertex (their difference is equal to the number of vertices that are on the same level with the given one and to the right of it).

4. This property follows from the consistency of edge orientations in the boundary of an elliptic region.

\begin{theorem}
\textbf{(Realization theorem)} If an ordered set of $n+1$ numbers is given, some of which have an overline and/or prime, and this set satisfies properties 1--4 of previous theorem, then it is the code of some 1-flow on the 2-disk.
\end{theorem}
\textbf{Proof.} We fix vertices on the boundary of the 2-disk. According to the given code, considering its successive numbers, we first draw neighboring loops at the given vertex in the amount equal to the first number. Then, inside the first loop, we draw loops in a quantity equal to the second number, inside the second loop - the third number, etc. Let's set the orientations of the drawn loops in accordance with the top lines. In the elliptical areas, we select the elliptical angles according to the primes. Thus, we have uniquely constructed a separatrix diagram that defines the flow.

$ $

To calculate the number of topologically non-equivalent flows, we find the number of such flows in each cell with a given orientation of the bottom side. Let the cell contain n top sides. If the orientations of all of them are opposite to the orientation of the bottom side, then the cell will be cyclic and each corner can be elliptical. In total we have n+1 such flow. If the orientation of at least one upper separatrix coincides with the orientation of the lower separatrix, then the cell is of polar type. The source angle can be selected in n ways. With the selected source angle, the sink angle can be located at any of the vertices on the path from the upper edges from the source angle to the end angle of the lower separatrix. Thus, there are 1+2+3+ ...+ n such configurations. Considering the n+1 configuration with a cyclic cell, we obtain that for a cell with a given orientation of the lower separatrix, 
$$1+2+\ldots +n+n+1=n(n+1)/2$$ 
different flow structures inside it are possible. Each flow structure in a cell defines a structure on the upper separatrices. Since these separatrices are lower for a cell of a higher level, then on this new cell we can assume that the orientation of the lower edge is given. Thus, the number of possible threads with a given graph (tree) is equal to the product of the number of threads in each cell.

Summing up the corresponding numbers in the last column of the table, we get the following result:

\begin{theorem}
 The number of topologically non-equivalent 1-flows on disk 2 without separatrices is 1, with one separatrix is 3, with two separatrices is 15, with three separatrices is  91, with four separatrices is 612, and with five separatrices is 4389.
\end{theorem}
\textbf{Remark.} Having carried out similar calculations for a flow with six separatrices, we obtain 31630 different structures. And for flows with seven separatrices we have 162900 different structures (see correspondent trees on Fig. \ref{trees6} and  Fig.\ref{trees7}).

\section*{Conclusion}
The results of the article prove the effectiveness of the proposed invariants (distinguishing graph and code) for classifying flows on a two-dimensional disk with one singular point. We hope that they can be used to construct similar invariants for flows with a large number of singular points on other surfaces as well.


\textsc{Taras Shevchenko National University of Kyiv}

Alexandr Prishlyak \ \ \ \textit{Email address:} \text{ prishlyak@knu.ua} \ \ \ \
\textit{ Orcid ID:} \text{0000-0002-7164-807X}

Serhii  Stas \ \ \ \  \ \ \ \ \ \ \ \ \  \ \textit{Email:} \text{ stasserhiy380@gmail.com} \  \ \ \ \ 
\textit{ Orcid ID:} \text{0009-0006-1241-7497}

\end{document}